\documentclass [12pt]  {article} 
\usepackage{amssymb}
\def \FF{{\mathbb{F}}}

\begin{document}

\begin{center}
{\Large {\bf On the exponent of the automorphism group of 
a compact Riemann surface}}\\
\bigskip
\bigskip
{\sc Andreas Schweizer\footnote{The author was supported by 
the National Research Foundation of Korea (NRF) grant funded
by the Korean government (MSIP) (ASARC, NRF-2007-0056093).}}\\
\bigskip
\bigskip
\it
To Ernst-Ulrich Gekeler in gratitude for his support 
and for everything that I learned from him.
\rm
\end{center}
\begin{abstract}
\noindent
Let $X$ be a compact Riemann surface of genus $g\geq 2$, and let
$Aut(X)$ be its group of automorphisms. We show that the exponent 
of $Aut(X)$ is bounded by $42(g-1)$. We also determine explicitly 
the infinitely many values of $g$ for which this bound is reached 
and the corresponding groups. Finally we discuss related questions 
for subgroups $G$ of $Aut(X)$ that are subject to additional 
conditions, for example being solvable.
\end{abstract}

\subsection*{1. Introduction}

\noindent
Throughout the paper $X$ will be a compact Riemann surface of 
genus $g\geq 2$. We write $Aut(X)$ for the full group of conformal 
automorphisms of $X$. 
\par
The order of a group $G$ is denoted by $|G|$, the cyclic group of 
order $n$ by $C_n$, and the neutral element in a group by $\iota$.
\\ \\
{\bf Theorem 1.1. (Hurwitz)} \it 
Let $X$ be a compact Riemann surface of genus $g\geq 2$. Then
$$|Aut(X)|\leq 84(g-1).$$
\rm

\noindent
Moreover, it is known that there are infinitely many values of $g$ for 
which the bound in Theorem 1.1 is reached. For proofs of all this see
for example [A, pp.46] or [Br, Theorem 3.17].
\par
A group $Aut(X)$ that reaches the bound $|Aut(X)|=84(g-1)$ is called 
a {\bf Hurwitz group}.
In Section 3 we will provide some more details on these groups.
\\ \\
After this one may of course ask how big a group $G$ that 
satisfies some additional properties can be if it acts as 
a group of automorphisms on a Riemann surface of genus $g$ 
(i.e. if $G\subseteq Aut(X)$). From the vast literature we 
select some results that are relevant for this paper.
\\ \\
{\bf Theorem 1.2.} \it 
Let $X$ be a compact Riemann surface of genus $g\geq 2$, and let
$G$ be a subgroup of $Aut(X)$.
\begin{itemize}
\item[(a)] If $G$ is solvable, then $|G|\leq 48(g-1)$. There are 
infinitely many values of $g$ for which this bound is reached.
\item[(b)] If $G$ is supersolvable, then for $g\geq 3$ we have
$|G|\leq 18(g-1)$. There are infinitely many values of $g$ for 
which this bound is reached. The biggest supersolvable group of 
automorphisms for genus $2$ has order $24$.
\item[(c)] If $G$ is nilpotent, then $|G|\leq 16(g-1)$. There are 
infinitely many values of $g$ for which this bound is reached.
\item[(d)] If $G$ is abelian, then $|G|\leq 4g+4$. For each 
$g\geq 2$ there are abelian groups of order $4g+4$ acting as 
automorphisms on a Riemann surface of genus $g$.
\end{itemize}
\rm

\noindent
{\bf Proof.} \rm 
(a) The bound results from the facts that Hurwitz groups are not 
solvable (see Corollary 3.2 (a) below) and that $48(g-1)$ is the 
next possible size of $Aut(X)$ [Br, Lemma 3.18]. Groups that reach 
the bound were constructed in [Ch] and [G1]. See also [G2] for 
improvements and minor corrections.
\par
(b) The papers [Z2] and [GMl] seem to have been written independently
and at almost the same time. In [Z2] the condition for $g$ to reach
the bound contains an error, which is pointed out in [GMl] and also 
corrected in [Z3].
\par
(c) [Z1, Theorems 1.8.4 and 2.1.2]
\par
(d) See [G1, p.271]. The paper [Ml] contains more precise information, 
namely on page 711 for each abelian group the minimal genus for which 
it can occur in $Aut(X)$.
\hfill$\Box$
\\ \\
Given $G\subseteq Aut(X)$, it is in general very difficult to decide 
whether $G$ equals $Aut(X)$ or is a proper subgroup. In the special case 
$|G|=48(g-1)$ we automatically have equality, as by [Br, Lemma 3.18] the 
only bigger order is $84(g-1)$, which is not a multiple of $48(g-1)$.
\par
A statement as in Theorem 1.2 for cyclic subgroups of $Aut(X)$ 
is of course equivalent to a statement about element orders. In 
[H, Theorem 6] for each $n$ the minimum genus for an automorphism 
of order $n$ is given. From this one can get the following classical 
result. Alternatively, see [G1, p.270].
\\ \\
{\bf Theorem 1.3. (Wiman)} \it 
Let $X$ be a compact Riemann surface of genus $g\geq 2$. Then the 
element orders of $Aut(X)$ are bounded by $4g+2$. 
\par
For each $g\geq 2$ there exists an $X$ of genus $g$ such that 
$Aut(X)$ contains elements of order $4g+2$.
\rm
\\ \\
The paper [N] classifies all Riemann surfaces with an automorphism
of order $\geq 3g$.
\\ \\
In this paper we investigate a problem that does not seem to have 
been treated yet in the literature, namely bounding the exponent
of $Aut(X)$ (or of a subgroup $G\subseteq Aut(X)$) in terms of the 
genus. 
\\

\subsection*{2. The exponent of a group}

\noindent 
The {\bf exponent} $exp(G)$ of a finite group $G$ is the least common 
multiple of all element orders. Equivalently, $exp(G)$ is the smallest 
positive integer $e$ such that $\sigma^e =\iota$ for all $\sigma\in G$. 
\\ \\
We leave the following three facts as easy exercises.
\\ \\
{\bf Lemma 2.1.} \it
$$exp(G)=\prod_p exp(G_p),$$
where the product is over all primes $p$ dividing $|G|$ and $G_p$
is a Sylow $p$-subgroup of $G$.
\rm
\\ \\
{\bf Lemma 2.2.} \it
The exponent of a finite $p$-group $P$ is the biggest element order.
In particular, $exp(P)=|P|$ if and only if $P$ is cyclic.
\rm
\\ \\
{\bf Corollary 2.3.} \it
\begin{itemize}
\item[(a)] $exp(G)=|G|$ if and only if all Sylow subgroups of 
$G$ are cyclic.
\item[(b)] $exp(G)=\frac{1}{2}|G|$ if and only if for all odd primes 
$p$ the Sylow $p$-subgroups of $G$ are cyclic and the (non-cyclic) 
Sylow $2$-subgroup has a cyclic subgroup of index $2$.
\end{itemize}
\rm

\noindent
Finite groups whose Sylow subgroups are all cyclic are called
{\bf Z-groups}, possibly from the German word zyklisch or perhaps
in honor of H. Zassenhaus, who completely described the structure
of these groups. For many purposes the form in [R, Theorem 10.1.10]
is better than the one in [Za, Satz 5].
\\ \\
{\bf Theorem 2.4. (Zassenhaus)} \it
A $Z$-group that is not cyclic can be written as a semidirect product
$$C_m\rtimes C_n$$
where $(m,n)=1$ and $m$ is odd. In particular, such a group is 
metacyclic and hence supersolvable.
\rm
\\ \\
The most important type of groups for our paper are the ones from
Corollary 2.3 (b), and the most important instance of such groups 
is the following.
\\ \\
{\bf Example 2.5.} 
Let $p$ be an odd prime. Then
$$|PSL_2(\FF_p)|=\frac{p^3 -p}{2}\ \ \ \hbox{\rm and} \ \ \ 
exp(PSL_2(\FF_p))=\frac{p^3 -p}{4}.$$

\noindent
Actually, finite groups $G$ with $exp(G)=\frac{1}{2}|G|$ have been
completely classified, the solvable ones in [Za, Satz 7] and the
nonsolvable ones in [S] and [W2]. See also the first page of [W1]
where it is explained why the different types discussed in [S] and
[W2] cover all possible cases. The final summary is
\\ \\
{\bf Theorem 2.6. (Suzuki, Wong)} [W2, Theorem 2] \it
Let $G$ be a nonsolvable finite group in which all Sylow subgroups 
of odd order are cyclic and a Sylow $2$-subgroup has a cyclic subgroup 
of index $2$. Then $G$ has a normal subgroup $G_1$ such that 
$[G:G_1]\leq 2$ and 
$$G_1 =L\times M,$$
where $L$ is isomorphic to $SL_2(\FF_p)$ or $PSL_2(\FF_p)$ for some
prime $p\geq 5$, and $M$ is a $Z$-group whose order is prime to that 
of $L$.
\rm
\\

\subsection*{3. Hurwitz groups}

In this section we collect the necessary details about Hurwitz groups.
See [C] for more background. 
\par
The first result comes from the fact that Hurwitz groups are exactly 
the non-trivial finite quotients of the triangle group 
$\Gamma(0;2,3,7)$. See for example [Br, Theorem 3.17]. Actually,
Section 3 of [Br] is a compact survey (with references) that covers
everything we need about triangle groups.
\\ \\
{\bf Theorem 3.1.} \it
A non-trivial, finite group is a Hurwitz group if and only if it can
be generated by two elements $\sigma$ and $\tau$ subject to
$$\sigma^2 =\tau^3 =(\sigma\tau)^7=\iota$$
and some other relations.
\rm
\\ \\
{\bf Corollary 3.2.} \it
Let $G$ be a Hurwitz group. Then
\begin{itemize}
\item[(a)] $G$ has no non-trivial abelian quotient group. So the 
commutator group $G'$ equals $G$, and $G$ is not solvable.
\item[(b)] Every non-trivial quotient group of $G$ has order divisible 
by $42$.
\item[(c)] Every non-trivial quotient of $G$ is again a Hurwitz group.
\end{itemize}
\rm

\noindent
{\bf Proof.} \rm 
We write $\tilde{\sigma}$ and $\tilde{\tau}$ for the images of
$\sigma$ and $\tau$ in the quotient.
\par
(a) If $\tilde{\sigma}$ and $\tilde{\tau}$ commute, then 
$(\tilde{\sigma}\tilde{\tau})^6=\tilde{\iota}$, 
and consequently $\tilde{\sigma}\tilde{\tau}=\tilde{\iota}$, 
$\tilde{\sigma}=\tilde{\iota}$, and $\tilde{\tau}=\tilde{\iota}$.
\par
(b) Similarly, if $42$ does not divide the order of the quotient
group, then one, and hence all, of $\tilde{\sigma}$, $\tilde{\tau}$, 
$\tilde{\sigma}\tilde{\tau}$ 
must equal $\tilde{\iota}$.
\par
(c) This is immediate from Theorem 3.1.
\hfill$\Box$
\\ \\
The question which finite groups are Hurwitz groups is far from 
being completely solved. Even for finite simple groups the answer 
is quite irregular (see [C]). We are interested in a special type 
of group.
\\ \\
{\bf Theorem 3.3. (Macbeath)} [Mb, Theorem 8] \it
The group $PSL_2(\FF_q)$ is a Hurwitz group if and only if
\begin{itemize}
\item[(i)] $q=7$,
\item[(ii)] $q=p$, a prime, with $p\equiv \pm 1\ (\hbox{\rm mod } 7)$,
\item[(iii)] $q=p^3$, where $p$ is a prime with $p\equiv \pm 2$ or
$\pm 3\ (\hbox{\rm mod } 7)$,
\end{itemize}
and for no other values of $q$.
\par
In cases (i) and (iii) there is only one Riemann surface on which
$G$ acts as a Hurwitz group. In case (ii) there are three Riemann
surfaces for each $G$.
\rm
\\ \\
In contrast, we point out the following easy result.
\\ \\
{\bf Theorem 3.4.} \it
$SL_2(\FF_p)$ is not a Hurwitz group.
\rm
\\ \\
{\bf Proof.} \rm 
Obviously the involution $\sigma$ (in Theorem 3.1) of a Hurwitz 
group cannot be central. But $SL_2(\FF_p)$ has exactly one involution, 
which thus is of course central. Actually, it is the negative of the 
unit matrix.
\hfill$\Box$
\\ \\
By Corollary 3.2 the exponent of a Hurwitz group has to be divisible 
by $42$. For use in later sections we refine this statement. To that
end we need the following group theoretic result.
\\ \\
{\bf Theorem 3.5.} \it 
Let $G$ be a non-abelian, simple group of order $2^a 3^b 5^c 7^d$ with 
abelian Sylow $2$-subgroup. Then $G$ must be among the groups 
$PSL_2(\FF_{p^n})$ with $p\in\{2,3,5,7\}$.
\rm
\\ \\
{\bf Proof.} \rm 
By [Wa, Theorem I] and the remarks immediately after it, a non-abelian 
finite simple group $G$ that has abelian Sylow $2$-subgroups and is not 
of type $PSL_2(\FF_{p^n})$ must either be the Janko group $J_1$ of order 
$175,560=2^3 \cdot 3\cdot 5\cdot 7\cdot 11\cdot 19$ or it must contain 
a subgroup $PSL_2(\FF_{3^{2n+1}})$ with $n>0$. 
\par
By an elementary number theoretic argument we show now that the order 
of such a group is always divisible by a prime $p>7$. Obviously,
$(3^{2n+1})^2 -1$ is congruent to $8$ (mod $16$), $-1$ (mod $3$), and
$3$ (mod $5$). So if it is not divisible by any prime $p>7$, it must
be of the form 
$$(3^{2n+1})^2 -1 = 8\cdot 7^m\ \ \ \hbox{\rm with}\ m>1.$$
Calculating modulo $9$ shows that $m$ must be divisible by $3$, 
say $m=3k$. And modulo $7$ we see that necessarily $n=3h+1$. Thus
$(3^{4h+2},2\cdot 7^k)$ is an integral solution of $x^3 -y^3 =1$.
But writing this as $(x-y)(x^2+xy+y^2)=1$ clearly shows the 
impossibility of such integral solutions.
\hfill$\Box$
\\ \\
Doubtlessly, the groups in Theorem 3.5 can be completely determined,
and most likely this is known. But the crude version above suffices 
for our goal, namely to prove
\\ \\
{\bf Theorem 3.6.} \it 
There are no Hurwitz groups of exponent $2\cdot 3\cdot 7^n$. In 
particular, there are no Hurwitz groups of exponent $42$.
\rm
\\ \\
{\bf Proof.} \rm 
Assume that $G$ is such a group. Let $N$ be a maximal normal subgroup.
By Corollary 3.2 (b) and (c), $G/N$ is a simple Hurwitz group, whose 
exponent is of the same form, possibly with a smaller $n$. 
\par
The Sylow $2$-subgroup of $G/N$ has exponent $2$; so in particular it 
must be abelian. Hence $G/N\cong PSL_2(\FF_{p^n})$ with $p\in\{2,3,5,7\}$
by Theorem 3.5. Now Theorem 3.3 leaves only the following four candidates 
for $G/N$, which however all fail: 
$PSL_2(\FF_7)$ has exponent $2^2\cdot 3\cdot 7$;
$PSL_2(\FF_8)$ has exponent $2\cdot 3^2 \cdot 7$;
$PSL_2(\FF_{27})$ has exponent $2\cdot 3\cdot 7\cdot 13$;
and $PSL_2(\FF_{125})$ has exponent $2\cdot 3^2 \cdot 5\cdot 7\cdot 31$.
\hfill$\Box$
\\ \\

\subsection*{4. The main result}

\noindent
{\bf Theorem 4.1.} (Main Theorem, first version) \it
Let $X$ be a compact Riemann surface of genus $g\geq 2$. Then
$$exp(Aut(X))\leq 42(g-1).$$
Equality holds if and only if $G=Aut(X)$ is a Hurwitz group
with $exp(G)=\frac{1}{2}|G|$.
\rm
\\ \\
{\bf Proof.} \rm 
If $exp(Aut(X))=|Aut(X)|$, then $Aut(X)$ is supersolvable by
Theorem 2.4. By Theorem 1.2 (b) therefore $exp(Aut(X))$ is 
significantly smaller than $42(g-1)$.
\par
In all other cases we have 
$exp(Aut(X))\leq \frac{1}{2}|Aut(X)|\leq \frac{84}{2}(g-1)$,
with equality if and only if $exp(Aut(X))=\frac{1}{2}|Aut(X)|$ 
and $|Aut(X)|=84(g-1)$.
\hfill$\Box$
\\ \\
{\bf Remark 4.2.} 
Theorem 4.1 already shows that there are infinitely many values of 
$g$ for which the bound $42(g-1)$ cannot be reached, simply because 
there are no Hurwitz groups for these $g$. Take for example $g=7^n +1$.
A group of order $84\cdot 7^n$ has a normal Sylow $7$-subgroup, so it
is solvable and thus cannot be a Hurwitz group. See [A, Chapter 5] for
more sequences of $g$ without Hurwitz groups.
\par
On the other hand, there are infinitely many values of $g$ for which 
at least one surface reaches $42(g-1)$. For example with the groups 
in Theorem 3.3 (i) and (ii). In the remainder of this section we want 
to show that these examples are the only ones.
\\ \\
{\bf Theorem 4.3.} \it
The only Hurwitz groups with $exp(G)=\frac{1}{2}|G|$ are the groups
$PSL_2(\FF_p)$ where $p=7$ or $p$ is a prime that is congruent to
$\pm 1$ modulo $7$.
\rm
\\ \\
{\bf Proof.} \rm 
As a complete classification of all Hurwitz groups is not known,
and would almost certainly be very complicated anyway, we start 
with the other condition. Let $G$ be a non-solvable group with 
$exp(G)=\frac{1}{2}|G|$. These are completely classified in
Theorem 2.6. 
\par 
If $G$ moreover is a Hurwitz group, we must have $G=G_1$ in that theorem 
by Corollary 3.2. Furthermore, the $Z$-group $M$ in that theorem must 
be trivial, for otherwise we could map from $G$ to $M$ and from there
to an abelian quotient, contradicting Corollary 3.2. So we are left with
the possibilities $G\cong PSL_2(\FF_p)$ or $G\cong SL_2(\FF_p)$. But the
second possibility is excluded by Theorem 3.4. Finally we apply
Theorem 3.3.
\hfill$\Box$
\\ \\
Combining Theorem 4.1 and Theorem 4.3 we obtain 
\\ \\
{\bf Theorem 4.4.} (Main Theorem, final version) \it
Let $X$ be a compact Riemann surface of genus $g\geq 2$. Then
$$exp(Aut(X))\leq 42(g-1).$$
This bound can be reached if and only if 
$$g=\frac{p^3 -p}{168}+1$$
where $p=7$ or $p$ is a prime that is congruent 
to $\pm 1$ modulo $7$. 
The only Riemann surface of genus $3$ that reaches the bound 
is the Klein quartic, whose automorphism group is isomorphic 
to $PSL_2(\FF_7)$. If $p\equiv\pm 1\ (\hbox{\rm mod } 7)$, 
there are $3$ non-isomorphic Riemann surfaces $X$ of genus 
$g=\frac{p^3 -p}{168}+1$ with $exp(Aut(X))=42(g-1)$. In every 
case
$$Aut(X)\cong PSL_2(\FF_p).$$
\rm
\\

\subsection*{5. Solvable groups}

In accordance with Theorem 1.2 we now try to find upper bounds on
$exp(G)$ for the $G\subseteq Aut(X)$ that are subject to additional
conditions.
\par
The following partial result, which might be interesting in its own 
right, will be used repeatedly.
\\ \\
{\bf Proposition 5.1.} \it
Let $X$ be a compact Riemann surface of genus $g\geq 2$, and let
$G\subseteq Aut(X)$ be a $Z$-group. Then
$$|G|<16(g-1).$$
\rm
\\ 
{\bf Proof.} \rm 
For $g=2$ there are four groups of order $\geq 16$ [Br, p.77], but
none of them has a cyclic Sylow $2$-subgroup.
\par
So we can suppose from now on that $g\geq 3$. Since $Z$-groups are of 
course metabelian, by [ChP] we have $|G|\leq 16(g-1)$ with two possible
exceptions, namely $|G|=48$ for $g=3$ and $|G|=80$ for $g=5$. But by
Theorem 1.2 (b) these two metabelian exceptions cannot be supersolvable,
and hence in particular not $Z$-groups.
\par
By [ChP] all metabelian groups $G$ of order $16(g-1)$ are quotients
of $\Gamma(0;2,4,8)$. If such a $G$ is a $Z$-group, then Theorem 2.4
implies that $G$, and hence also $\Gamma(0;2,4,8)$ has a cyclic quotient
of order $16$, which is clearly impossible.
\hfill$\Box$
\\ \\
{\bf Remark 5.2.} 
We don't know what could be a sharp bound in Proposition 5.1. In any
case there are infinitely many $Z$-groups of order $10(g-1)$. Namely,
by [BJ, Theorem 1] for every big enough prime $p$ with 
$p\equiv 1\ (mod\ 5)$ there exists a Riemann surface $X$ of genus 
$p+1$ such that $Aut(X)$ contains $G\cong C_p\rtimes C_{10}$.
\\ \\
{\bf Proposition 5.3.} \it
Let $X$ be a compact Riemann surface of genus $g\geq 3$, and 
let $G$ be a solvable subgroup of $Aut(X)$. Then
$$exp(G)\leq 16(g-1).$$
\rm
\\ 
{\bf Proof.} \rm 
If $exp(G)<\frac{1}{2}|G|$, Theorem 1.2 (a) implies 
$exp(G)\leq \frac{1}{3}|G|\leq16(g-1)$. 
And if $exp(G)=|G|$, we even have $exp(G)<16(g-1)$ 
by Proposition 5.1.
\par
Solvable groups $G$ with $exp(G)=\frac{1}{2}|G|$ have been 
completely classified in [Za, Satz 7]. We only need the 
following key fact from the proof: If $exp(G)=\frac{1}{2}|G|$, 
then $G$ contains a normal $Z$-group $G_1$ such that $G/G_1$ 
is isomorphic to $C_2$ or $A_4$ or $S_4$.
\par
If $[G:G_1]=2$, we have $exp(G)=|G_1|<16(g-1)$ by Proposition 5.1.
In the remaining cases we have to show that $|G|\leq 32(g-1)$. By
[Br, Lemma 3.18] there are only three possible orders of solvable 
groups bigger than $32(g-1)$, namely $48(g-1)$, $40(g-1)$ and 
$36(g-1)$. Correspondingly, we have to show that $G$ cannot be 
a quotient of $\Gamma(0;2,3,8)$, $\Gamma(0;2,4,5)$ or 
$\Gamma(0;2,3,9)$.
\par
Obviously, the only finite quotient of $\Gamma(0;2,4,5)$ of 
order prime to $5$ is $C_2$. So $\Gamma(0;2,4,5)$ cannot have 
a quotient $G$ that has a quotient $A_4$ or $S_4$. 
Likewise, $\Gamma(0;2,3,8)$ has no quotient $C_3$ and hence 
no quotient $G$ with $G/G_1\cong A_4$, whereas 
$\Gamma(0;2,3,9)$ has no quotient $C_2$ and hence no quotient 
$G$ with $G/G_1\cong S_4$, 
\par
For the remaining two cases we use that since $G_1$ is 
supersolvable, the elements in $G_1$ of odd order form 
a characteristic subgroup $U$ of $G_1$ [R, Theorem 5.4.9].
The normality of $G_1$ in $G$ implies that $U$ is normal 
in $G$. Furthermore, if $U$ is non-trivial and $p_1$ is the 
smallest prime divisor of $|U|$, by Zappa's Theorem 
[R, Theorem 5.4.8] $U$ has a normal subgroup $M$ of index
$p_1$. Since the Sylow $p_1$-subgroups of $U$ are cyclic,
$M$ is even characteristic in $U$, and hence normal in $G$.
\par
If $G/G_1\cong S_4$ we obtain a chain of normal subgroups
$$G/M\triangleright N\triangleright V\triangleright G_1/M\triangleright
U/M\triangleright I$$
with factors $C_2$, $C_3$, $C_2\times C_2$, $C_{2^e}$, $C_{p_1}$.
Now assume moreover that $G$ is a quotient of $\Gamma=\Gamma(0;2,3,8)$.
By [Br, Example 3.8] we have $\Gamma'\cong\Gamma(0;3,3,4)$ and
$\Gamma''\cong\Gamma(0;4,4,4)$ with $\Gamma/\Gamma'\cong C_2$ and
$\Gamma'/\Gamma''\cong C_3$. This implies that $N$ is a quotient 
of $\Gamma'$ and $N'=V$ and $V$ is a quotient of $\Gamma''$.
In particular, $V$ cannot have a quotient $C_{p_1}$.
\par
On the other hand, $U/M$ is a normal subgroup of $N$. Let $C$ be 
its centralizer in $N$. As $N/C$ can be embedded into the automorphism 
group of $U/M$, which is cyclic, we see $C\supseteq N'=V$. So $U/M$ is 
central in $V$. This means that $V$ is a direct product of its Sylow 
$2$-subgroup and $U/M$. In particular, $V$ has a quotient $C_{p_1}$.
\par
The resolution of this contradiction is that $U$ must be trivial.
Consequently $48(g-1)=|G|=3\cdot 2^{e+3}$. Since $exp(G)=24(g-1)$, 
in that case $G$ must contain an element of order $8(g-1)$. By 
Theorem 1.3 this is only possible if $8(g-1)\leq 4g+2$, i.e., 
if $g\leq 2$.
\par
Similarly, if $G$ is a quotient of $\Gamma=\Gamma(0;2,3,9)$ with
$G/G_1\cong A_4$, we obtain that $G'$ lies between $G$ and $G_1$ 
with $G/G'\cong C_3$ and $G'$ is a quotient of the commutator group
$\Gamma'\cong\Gamma(0;2,2,2,3)$. Moreover, since $36$ divides $|G|$,
we have $U/M\cong C_3$, and as above this group is central in $G'/M$, 
leading to the contradiction that $\Gamma'$ should have a quotient 
$C_3$.
\hfill$\Box$
\\ \\
We don't know whether for $g>2$ the bound in Proposition 5.3 can 
be reached, and if yes whether infinitely often.
\par
For genus $2$ we mention that the Bolza surface $y^2 = x^5 -x$ has 
automorphism group $GL_2(\FF_3)$ of order $48$ and exponent $24$. 
\\ \\
{\bf Remark 5.4.} 
By [Z2, Theorem 4.1] or [GMl, Lemma 4.1] the supersolvable group of 
order $24$ for $g=2$ has exponent $12$.
\par
If $g\geq 3$ and $G$ is supersolvable but not a $Z$-group, then from
Theorem 1.2 (b) we get $exp(G)\leq \frac{1}{2}|G|\leq 9(g-1)$, which
is smaller than the examples mentioned in Remark 5.2. This shows that 
bounding $exp(G)$ for supersolvable groups $G$ amounts to the same as 
bounding $|G|$ for $Z$-groups $G$.
\\ \\
A finite nilpotent group $G$ is the direct product of its Sylow 
subgroups. Its exponent therefore is the biggest element order 
(see Lemmas 2.1 and 2.2). So by Theorem 1.3 we have 
$exp(G)\leq 4g+2$. By [N, Theorem 1] there is a unique surface 
of genus $g$ that has an automorphism of order $4g+2$. Putting 
all together we obtain
\\ \\
{\bf Theorem 5.5.} \it
Let $X$ be a compact Riemann surface of genus $g\geq 2$, and let
$G$ be a nilpotent subgroup of $Aut(X)$. Then
$$exp(G)\leq 4g+2.$$
For every $g\geq 2$ there exists, up to isomorphism, exactly one
Riemann surface of genus $g$ that realizes this bound, namely 
$$y^2 =x^{2g+1}-1,$$
which has $Aut(X)\cong C_{4g+2}$.
\rm
\\ \\
Obviously the same result holds for abelian subgroups $G$ of 
$Aut(X)$.
\\

\subsection*{6. On $|G|/exp(G)$}

Finally, we investigate the case when $exp(G)$ is as small as
possible compared to $|G|$.
\\ \\
{\bf Theorem 6.1.} \it 
Let $X$ be a compact Riemann surface of genus $g\geq 2$, and let
$G$ be a subgroup of $Aut(X)$. Then 
$$\frac{|G|}{exp(G)}\ \ \hbox{\it divides}\ \ 2(g-1).$$
\rm
\\ 
{\bf Proof.} \rm 
This is a well-known consequence of the Hurwitz formula 
$$2g-2=|G|(2h-2)+\sum_{i=1}^r \frac{|G|}{|S_i|}(|S_i|-1)$$
for the covering $X\to X/G$. Here $P_1 ,\ldots , P_r$ are the branch 
points on the genus $h$ Riemann surface $X/G$, and $S_i$ is the 
stabilizer of a point on $X$ above $P_i$. Since $S_i$ is always cyclic, 
$\frac{|G|}{exp(G)}$ divides $\frac{|G|}{|S_i|}$.
\hfill$\Box$
\\ \\
{\bf Proposition 6.2.} \it 
Let $X$ be a compact Riemann surface of genus $g\geq 2$, and let
$G$ be a subgroup of $Aut(X)$. If $|G|/exp(G)=2(g-1)$, then $G$
must be solvable and $exp(G)\leq 24$.
\rm
\\ \\
{\bf Proof.} \rm 
If $|G|/exp(G)=2(g-1)$, then Theorem 1.1 implies $exp(G)\leq 42$, 
with equality if and only if $G$ is a Hurwitz group of exponent 
$42$. But by Theorem 3.6 such groups do not exist. 
\par
The second biggest possible size of $G$ is $|G|=48(g-1)$ [Br, Lemma 3.18]. 
This shows $exp(G)\leq 24$. Consequently, $exp(G)$ cannot have more than 
$2$ different prime divisors. By Burnside's $p^m q^n$-Theorem (see for 
example [R, Theorem 8.5.3]) this implies that $G$ is solvable.
\hfill$\Box$
\\ \\
If $|G|/exp(G)=2(g-1)$, then $exp(G)$ must of course be even and
divisible by all primes that divide $g-1$. 
\par
Moreover, the cases $exp(G)=22$, $16$ or $14$ cannot occur in 
Proposition 6.2, because there are no $G$ of order $44(g-1)$,
$32(g-1)$ or $28(g-1)$ ([Br, Lemma 3.18]). 
On the other hand, the Bolza surface from Remark 5.2 shows that 
$exp(G)=24$ can occur, at least for $g=2$. We don't know whether
it can occur for $g>2$. But we have the following general 
finiteness result.
\\ \\
{\bf Theorem 6.3.} \it 
There are only finitely many groups $G$ that reach the bound 
$|G|/exp(G)=2(g-1)$ in Theorem 6.1.
\rm
\\ \\
{\bf Proof.} \rm 
Let $|G|/exp(G)=2(g-1)$. Dividing the Hurwitz formula in the proof 
of Theorem 6.1 by $|G|/exp(G)$, we obtain
$$1=exp(G)(2h-2)+\sum_{i=1}^r \frac{exp(G)}{|S_i|}(|S_i|-1).$$
Fix one of the remaining exponents $24$, $20$, $18$, $12$, $10$, 
$8$, $6$, $4$, $2$. Then, since $|S_i|$ divides $exp(G)$, there are 
only finitely many values $h$ and $r$ for which this equation has 
a solution. More precisely, we must have $h=0$ and $r\leq 5$, as 
$h=1$ would imply $r=1$, $|S_1|=2$, and $exp(G)=2$, which is not 
possible. By the theory of Fuchsian groups, $G$ is a quotient of
a group that is generated by $2h+r-1$ elements. But by the affirmative 
solution to the restricted Burnside problem [Ze] there are only 
finitely many finite groups with a given number of generators and 
a given exponent.
\hfill$\Box$ 
\\ \\
Finally, the case $exp(G)=2$ can be completely settled. All groups 
of exponent $2$ are abelian; so we do this in slightly more generality.
\\ \\
{\bf Theorem 6.4.} \it 
There are only five abelian groups $G$ that reach the bound 
$|G|/exp(G)=2(g-1)$ in Theorem 6.1, namely
\begin{itemize}
\item $C_2 \times C_2$ and $C_6 \times C_2$ for $g=2$;\
\item $C_2 \times C_2 \times C_2$ and $C_4 \times C_4$ for $g=3$;
\item $C_2 \times C_2 \times C_2 \times C_2$ for $g=5$.
\end{itemize}
\rm

\noindent
{\bf Proof.} \rm 
If $G\subseteq Aut(X)$ is abelian, then $|G|\leq 4g+4$ by 
Theorem 1.2 (d). If moreover $|G|/exp(G)=2(g-1)$, this leaves 
only the possibilities $exp(G)\in\{2,4,6\}$ for $g=2$, 
$exp(G)\in\{2,4\}$ for $g=3$, and $exp(G)=2$ for $g>3$. 
So besides the groups listed in the theorem, the possible 
candidates are $C_4 \times C_2$ for $g=2$, $C_4 \times C_2 \times C_2$
for $g=3$, and $(C_2)^r$, $r\geq 5$ for $g=2^{r-2}+1$.
But [Ml, Theorem 4] shows that the minimum genus for the latter
three types is bigger.
\hfill$\Box$ 
\\ \\ \\
{\bf Acknowledgements.} I am indebted to the anonymous referee for 
her/his careful reading of the manuscript and for many constructive
comments. These have not only improved the presentation, but have 
also substantially strengthened some of the results.
Specifically, the trick with the centralizer of $U/M$ in the proof
of Proposition 5.3 is due to the referee. This has improved the bound 
from originally $24(g-1)$ to $16(g-1)$. 
Equally due to the referee is the observation that in Theorem 6.3 the
number of generators of $G$ can be bounded for every exponent, which
is crucial for getting finiteness overall.
\\

\subsection*{\hspace*{10.5em} References}
\begin{itemize}

\item[{[A]}] R.~Accola: \it Topics in the Theory of Riemann Surfaces,
\rm Springer Lecture Notes in Mathematics 1595, 
Berlin-Heidelberg-New York, 1994

\item[{[BJ]}] M.~Belolipetsky and G.~A.~Jones: Automorphism 
groups of Riemann surfaces of genus $p+1$, where $p$ is prime,
\it Glasgow Math. J. \bf 47 \rm (2005), 379-393

\item[{[Br]}] T.~Breuer: \it Characters and Automorphism Groups 
of Compact Riemann Surfaces, \rm LMS Lecture Notes 280, Cambridge 
University Press, Cambridge, 2000

\item[{[Ch]}] B.~P.~Chetiya: On genuses of compact Riemann 
surfaces admitting solvable automorphism groups, 
\it Indian J. Pure Appl. Math. \bf 12 \rm (1981), 1312-1318 

\item[{[ChP]}] B.~P.~Chetiya and K.~Patra: On metabelian groups of 
automorphisms of compact Riemann surfaces, \it J. London Math. Soc.
\bf 33 \rm (1986), 467-472

\item[{[C]}] M.~Conder: Hurwitz groups: a brief survey, \it 
Bull. Amer. Math. Soc. (N.S.) \bf 23 \rm (1990), 359-370 

\item[{[G1]}] G.~Gromadzki: Maximal groups of automorphisms of
compact Riemann surfaces in various classes of finite groups,
\it Rev. Real Acad. Cienc. Exact. Fis. Natur. Madrid \bf 82 no. 2
\rm (1988), 267-275

\item[{[G2]}] G.~Gromadzki: On soluble groups of automorphism of 
Riemann surfaces, \it Canad. Math. Bull. \bf 34 \rm (1991), 67-73

\item[{[GMl]}] G.~Gromadzki and C.~Maclachlan: Supersoluble
groups of automorphisms of compact Riemann surfaces, 
\it Glasgow Math. J. \bf 31 \rm (1989), 321-327

\item[{[H]}] W.~J.~Harvey: Cyclic groups of automorphisms 
of a compact Rieman surface, \it Quart. J. Math. Oxford (2), 
\bf 17 \rm (1966), 86-97

\item[{[Mb]}] A.~M.~Macbeath: Generators of the linear fractional groups,
in: Number Theory (Houston 1967), \it Proc. Sympos. Pure Math. \bf 12,
\rm American Mathematical Society, Providence (1969), 14-32 

\item[{[Ml]}] C.~Maclachlan: Abelian groups of automorphisms of compact 
Riemann surfaces, \it Proc. London Math. Soc. (3) \bf 15 \rm (1965),
699-712 

\item[{[N]}] K.~Nakagawa: On the orders of automorphisms of a closed
Riemann surface, \it Pacific J. Math. \bf 115 no. 2 \rm (1984), 435-443 

\item[{[R]}] D.~J.~S.~Robinson: \it A Course in the Theory of Groups,
\rm Springer GTM 80, New York - Berlin, 1982

\item[{[S]}] M.~Suzuki: On finite groups with cyclic Sylow subgroups 
for all odd primes, \it Amer. J. Math. \bf 77 \rm (1955), 657-691 

\item[{[Wa]}] J.~H.~Walter: The characterization of finite groups
with abelian Sylow $2$-subgroups, \it Ann. of Math. (2) \bf 89 \rm 
(1969), 405-514 

\item[{[W1]}] W.~J.~Wong: On finite groups whose $2$-Sylow subgroups
have cyclic subgroups of index $2$, \it J. Austral. Math. Soc. \bf 4
\rm (1964), 90-112 

\item[{[W2]}] W.~J.~Wong: On finite groups with semi-dihedral Sylow 
$2$-subgroups, \it J. Algebra \bf 4 \rm( 1966), 52-63 

\item[{[Za]}] H.~Zassenhaus: \"Uber endliche Fastk\"orper, \it
Abh. Math. Sem. Univ. Hamburg \bf 11 \rm (1935), 187-220 

\item[{[Ze]}] E.~I.~Zelmanov: A solution of the restricted Burnside 
problem for $2$-groups, \it Math. USSR Sbornik \bf 72 \rm (1991), 
543-565 

\item[{[Z1]}] R.~Zomorrodian: Nilpotent automorphism groups of 
Riemann surfaces, \it Trans. Amer. Math. Soc. \bf 288 no. 1 
\rm (1985), 241-255

\item[{[Z2]}] R.~Zomorrodian: Bounds for the order of supersoluble
automorphism groups of Riemann surfaces, \it Proc. Amer. Math. Soc.
\bf 108 no. 3 \rm (1990), 587-600

\item[{[Z3]}] R.~Zomorrodian: On a theorem of supersoluble 
automorphism groups, \it Proc. Amer. Math. Soc.
\bf 131 no. 9 \rm (2003), 2711-2713

\end{itemize}

\noindent
{\sc Andreas Schweizer}\\
Department of Mathematics,\\
Korea Advanced Institute of Science and Technology (KAIST),\\ 
Daejeon 305-701,\\
South Korea\\
e-mail: {\tt schweizer@kaist.ac.kr}

\end{document}